# Market Integration of HVDC Lines:
# Cost Savings from Loss Allocation and Redispatching

A. TOSATTO[1]*, M. DIJOKAS[1], D. OBRADOVIC[2], T. WECKESSER[3], R. ERIKSSON[4], J. JOSEFSSON[5], A. KRONTIRIS[6], M. GHANDHARI[2], J. ØSTERGAARD[1], S. CHATZIVASILEIADIS[1]

[1]Technical University of Denmark (DTU), [2]KTH Royal Institute of Technology, [3]Dansk Energi, [4]Svenska kraftnät, [5]Hitachi ABB, [6]Darmstadt University of Applied Sciences
Denmark, Sweden, Germany

## SUMMARY

In the last decades, over 25'000 km of High-Voltage Direct-Current (HVDC) lines have been gradually integrated to the existing pan-European HVAC system. From a market point of view, HVDC interconnectors facilitate the exchange of energy and ancillary services between countries. In the Nordic region, many interconnectors are formed by HVDC links, as Scandinavia, Continental Europe and the Baltic region are non-synchronous AC systems. In this regard, this paper presents two cost benefit analyses on the utilization of HVDC interconnectors in the Nordic countries: in the first we investigate the utilization of HVDC interconnectors for reserve procurement and, in the second, we assess the implementation of implicit grid losses on HVDC interconnectors in the day-ahead market.

The first analysis is motivated by real events in 2018 where the inertia of the Nordic system dropped below a critical level and the most critical generating unit, a nuclear power plant in Sweden, was redispatched to guarantee the security of the system. In order to guarantee system security while reducing the costs of preventive actions, in summer 2020 new frequency products were introduced in the Nordic system: the Fast Frequency Reserves (FFR). HVDC lines, however, can perform similar tasks at lower costs. In our analysis, we are, thus, investigating the cost savings of using HVDC lines for frequency support using their Emergency Power Control (EPC) functionality, instead of redispatching or FFR.

The second analysis is based on the proposition of Nordic Transmission System Operators (TSOs) to introduce linear HVDC loss factors in the market clearing. With our analysis, we show that linear loss factors can unfairly penalize one HVDC line over the other, and this can reduce social benefits and jeopardize revenues of merchant HVDC lines. In this regard, we propose piecewise-linear loss factors: a simple-to-implement but highly-effective solution. Moreover, we demonstrate how the introduction of HVDC loss factors is a partial solution, since it disproportionally increases the AC losses. Our results show that the additional inclusion of AC loss factors can eliminate this problem.

* antosat@elektro.dtu.dk

## I. INTRODUCTION

The world's first commercial High-Voltage Direct-Current (HVDC) link was delivered by ABB in the 1950's. Since then, HVDC has become a common tool in the design of transmission grids, especially when technical limitations of AC transmission come into play. Indeed, the transmission of power in the DC form presents several benefits, both from technical and economical points of view [1], [2]. First, beyond a certain distance, an HVDC line has lower power losses than an HVAC of the same capacity. Second, with HVDC no reactive compensation is needed, resulting in no length limitation for submarine or underground power cables. Moreover, HVDC enables the connection of non-synchronous areas, allowing both inter-area and cross-continental long-distance power flows.

Evidence of the economic value of HVDC can be found by looking at the evolution of prices in different electricity markets. Installing transmission capacity means allowing for power exchanges between low- and high-price areas, reducing price differences and increasing social welfare. For example, Storebælt, the connection between the two bidding zones in Denmark (DK1 and DK2), has decreased electricity prices in Eastern Denmark by 2 €/MWh in average since 2010, resulting in 20-25 million euro savings per year for Danish consumers [3]. Another example is NordBalt, the link between Sweden and Lithuania: electricity prices in Lithuania dropped by 30% when the link was operated for the first time in 2016, followed by an average decrease of 5 €/MWh compared to before 2016 [4].

Furthermore, HVDC is attracting increasing attention because of the full controllability of power flows; depending on the technology of the converter stations, both active and reactive power flows can be controlled. This property can help improve the performance of AC power systems by means of additional control facilities, such as in the provision of frequency support.

Frequency stability is becoming a concern for many Transmission System Operators (TSOs) because of the gradual decrease of system inertia. This is mainly caused by the replacement of conventional synchronous generators with inverter-based non-synchronous units. Therefore, power systems are becoming more sensitive to power disturbances and TSOs have to bear additional costs for system security [5]. Until May 2020, the standard procedure to ensure N-1 security during low inertia periods was the reduction of the dimensioning incident (DI) – the largest disturbance. In Regional Group Nordic (RG Nordic), this is the loss of the most critical generating unit, a 1450 MW nuclear power plant in Sweden – Oskarshamn 3 (O3). For example, during Summer 2018, the power output of O3 was reduced by 100 MW three times [6]. To avoid such costly redispatching actions, since summer 2020, Fast Frequency Reserves (FFR) have been introduced as an ancillary services market product, which Nordic TSOs procure to complement the response of Frequency Containment Reserves for Disturbances (FCR-D) during low inertia periods. The costs of these mitigation strategies amount respectively to 0.38 million euros in 2018 (DI reduction) and 3.33 million Euros in 2020 (FFR). This trend of increasing cost for ancillary services has been observed in other countries, too, and is expected to continue with increasing shares of renewable generation [7], calling for a reassessment of whether there exist more cost-efficient options to guarantee safe operation while avoiding expensive remedial actions.

According to [8], the control scheme of all HVDC converters must be capable to operate in frequency sensitive mode, i.e. the transmitted power is adjusted in response to a frequency deviation. For this reason, an HVDC link connecting asynchronous areas can be used as a vehicle for sharing frequency reserves among asynchronous areas: to limit the instantaneous frequency deviation (IFD) in case of disturbance, the necessary active power can be imported from the neighboring system using the Emergency Power Control (EPC) functionality [9]. Given the high number of interconnections formed by HVDC links between RG Nordic and the neighboring regional groups, this corrective action could represent a valid alternative to expensive preventive redispatching or the procurement of expensive frequency reserves.

On the one hand, HVDC interconnectors are of great value for society as they facilitate the exchange of energy and ancillary services between countries. On the other hand, HVDC operation comes with a cost for TSOs as HVDC interconnectors produce a non-negligible amount of losses that is currently not considered in the market clearing. During periods of zero price difference between neighboring bidding zones, due to equal zonal prices, the cost of HVDC losses is transferred to local Transmission System Operators (TSOs) who must procure sufficient power to cover these losses. The problem is especially pronounced in transit countries, as in the case of Denmark.

Power losses are generated both by AC and DC interconnectors; however, because of their properties, HVDC lines are often significantly longer than AC lines and thus the operation of such lines leads to a considerable amount of losses. In 2017, the total losses by all the HVDC links in the Nordic region was equal 1.14 TWh [4]. Compared to the total amount of losses in the Nordic AC system, which is about 10 TWh, this value accounts for 10% [10]. However, given that all HVDC lines connect control areas operated by different TSOs, it is often unclear who should pay for these losses.

Recently, Nordic TSOs have proposed the introduction of HVDC loss factors (also called "implicit grid loss"] to implicitly account for losses when the market is cleared [11]. The introduction of loss factors will force a price difference between the two connected bidding zones that is equal to the marginal cost of losses. This will have two advantages: first, HVDC losses are no longer needed to be purchased by TSOs in the day-ahead market but are directly paid by the market participants who create them and, second, losses are implicitly minimized, resulting in cost savings for TSOs and the society. The proposed loss factors are linear approximations of the HVDC system losses. The following questions arise: are linear loss factors a good representation of HVDC losses? Is the introduction of loss factors for only HVDC interconnectors the best possible action?

In this regard, this paper presents a cost benefit analysis on:

- The utilization of HVDC interconnectors for frequency support, using the Emergency Power Control functionality to fulfil the N-1 security criterion.
- The implementation of implicit grid losses on HVDC interconnectors in the day-ahead market.

For the first analysis, we start by investigating what is the cost of frequency balancing using HVDC in the form of EPC, and then compare this alternative to DI reduction and FFR. The analysis is carried out for two scenarios for the year 2025, using inertia forecasts from Nordic TSOs. From these, the volume of EPC, FFR and redispatch is estimated, and the costs of these three mitigation strategies are compared.

For the second analysis, we compare the results of different simulations where the day-ahead market is cleared for each hour of the year (8760 instances) using data from 2017. Each simulation is carried out with different combinations of AC and HVDC loss factors, using different linearization techniques.

The rest of the paper is organized as follows. Section II provides background information about the Nordic power system and electricity market. Section III presents the results of the cost benefit analysis on the utilization of HVDC lines for frequency support and Section IV presents the analyses on the introduction of loss factors in the Nordics. Section V gathers conclusions and final remarks.

## II. NORDIC POWER SYSTEM

The Nordic transmission network is divided into two asynchronous Regional Groups (RGs): Western Denmark is connected to Continental Europe (UCTE) and, thus, it is operated at a different frequency from the rest of the Nordic countries. A schematic representation of the transmission network is depicted in Figure 1 (left).

Western Denmark is connected to Germany through different AC lines, along a corridor which is usually referred to as east coast corridor. Three HVDC links (Skagerrak, Kontiskan and Storebælt) connect Western Denmark to Norway, Sweden and Eastern Denmark. Recently, COBRAcable HVDC link has become operational, allowing power exchanges between Western Denmark and the Netherlands.

RG Nordic is connected to RG Continental Europe through five additional HVDC links: NorNed (Norway-Netherlands), Kontek (Eastern Denmark-Germany), Baltic cable (Sweden-Germany) and SwePol (Sweden-Poland). In addition, Kriegers Flak (Combined Grid Solution) provides connection between Eastern Denmark and Germany (AC cable and back-to-back HVDC converter), integrating offshore wind farms along its path. Finally, three other HVDC links connect RG Nordic to RG Baltic: NordBalt (Sweden-Lithuania), Estlink (Finland-Estonia) and Vyborg HVDC (Finland-Russia).

The generation mix in the Nordic countries can be found in [12]. Almost half of the generation in Denmark comes from wind farms, while the remaining is mainly fossil fuel based (natural gas and coal). In Norway, more than 90% of electricity is produced by hydro power plants. Hydro power plants contribute to half of the generation in Sweden as well, the remaining capacity is divided between nuclear power plants, wind farms and oil-based thermal units. In Finland the generation mix is more heterogeneous; half of the Finnish electrical energy is produced by nuclear power plants and coal-based thermal units.

As for the rest of Europe, the system is operated at 50 Hz with a standard range of ±100 mHz; Frequency Containment Reserves for Normal operation (FCR-N) are deployed to keep frequency within the normal band [13]. When frequency drops below 49.9 Hz, FCR for Disturbance (FCR-D) are activated to mitigate the impact of the disturbance and stabilize the frequency, while Frequency Restoration Reserves (FRR) are used to restore the frequency back to the nominal value. The maximum acceptable Instantaneous Frequency Deviation (IFD) is 1000 mHz and, in case frequency drops below 48.8 Hz, loads are shed to avoid total system blackout [14].

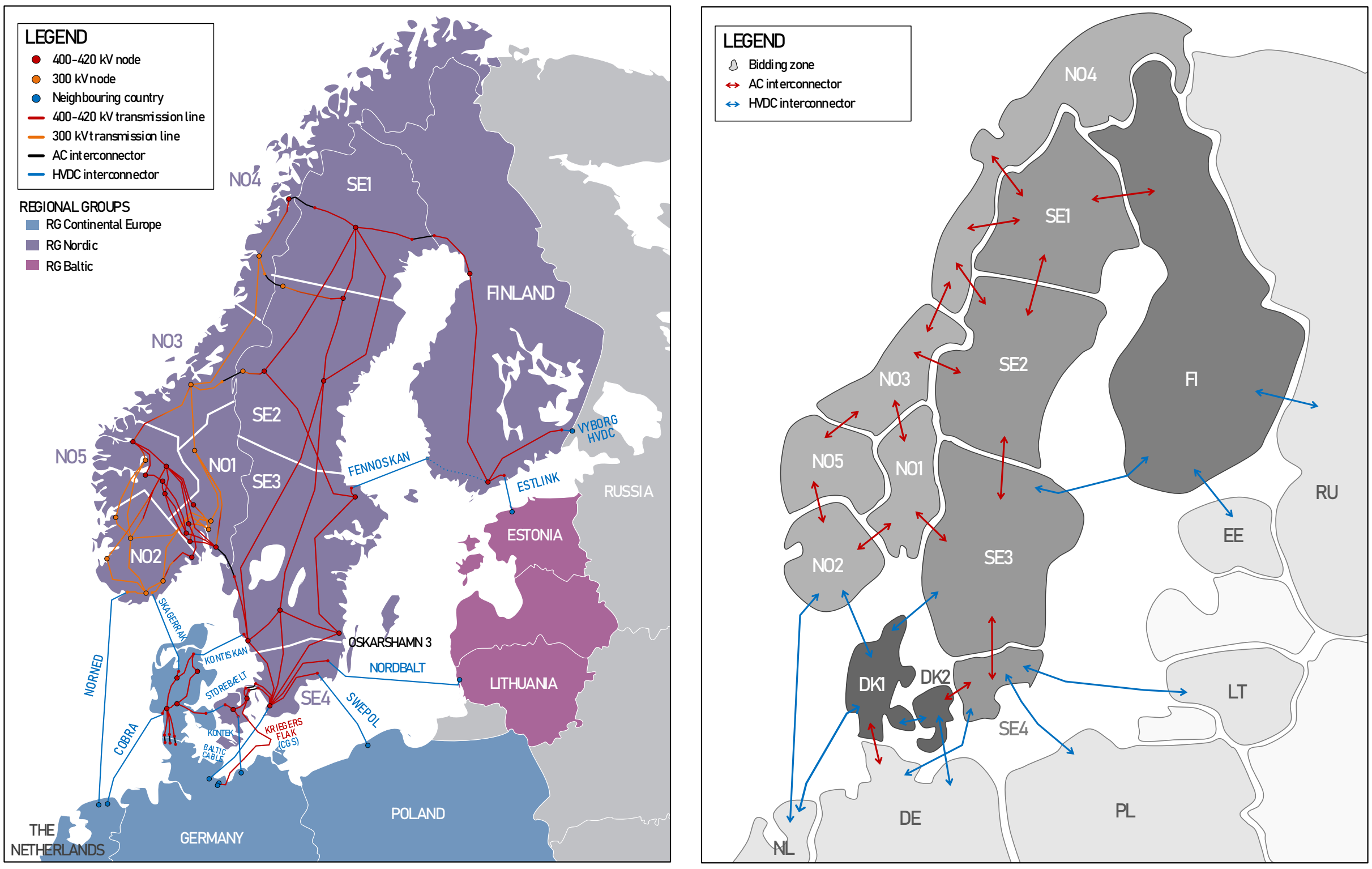

Figure 1: Nordic power grid and regional groups (left), Nordic market model (right).

The methodology for calculating the FCR-D requirement consists in a probabilistic approach which aims at reducing the probability of insufficient reserves, based on different generation, load and inertia patterns [15]. The considered dimensioning incidents are the loss of critical components of the system, such as large generators, demand facilities and transmission lines. Currently, the dimensioning incident in RG Nordic is the loss of Oskarshamn 3, a 1450 MW nuclear power plant in Sweden (located in the bidding zone SE3) [16].

Finally, in the Nordic region, as for the rest of Europe, a zonal-pricing scheme is applied. This means that the system is split into several bidding zones and the intra-zonal network is not included in the market model. When the market is cleared, a single price per zone is defined. In case of congestion, price differences arise only among zones [17]. The current day-ahead market coupling is based on Available Transfer Capacity (ATC). In the day-ahead time frame, TSOs calculate ATCs based on the network situation and communicate them to the market operator. These values are used as bounds for inter-zonal power transfers in the spot-market. When the power exchanges are defined, TSOs manage the physical flows to guarantee these transactions and, if necessary, counter-trade at their own cost [18]. Figure 1 (right) shows the different bidding zones in the Nordic area and the equivalent interconnectors.

## III. SHARING RESERVES THROUGH HVDC INTERCONNECTORS

In 2018, the power output of O3 has been reduced three times, due to system inertia dropping below the acceptable limit. Such low inertia periods are considered extraordinary events where the security of the system is in danger, thus Svenska kraftnät can communicate the limitation on O3 at any market stage. When this happens, the producer (in this case Oskarshamn Group – OKG) should receive market compensation for the costs associated with the power limitation. First, by decreasing its power output, the producer incurs opportunity costs that are equal to what they would have received for producing an amount of power equal to the power reduction. Second, by moving away from the nominal power output, extra costs are incurred due to lower efficiency (as a rule of thumb, for nuclear power plants, one can say that half of the fuel which is not used during the power reduction is lost and cannot be used later on) [9]. Third, the decrease of power production of nuclear power plants results in a temperature transient, inducing a cumulative aging of the unit and increasing the risk of failure [9]. All these factors are taken into consideration by the Nordic TSOs and OKG receives a financial compensation, as stated on a bilateral agreement between Svk and OKG. For each event of 2018:

- Oskarshamn 3 was compensated for the opportunity cost of not producing 100 MW (the compensation was equal to 49 SEK/MWh - approx. 4.64 €/MWh).
- Oskarshamn 3 was compensated for reduced efficiency and other costs associated with the power limitation (fixed amount equal to 50'000 SEK - approx. 4'740 €).

The downregulation was announced by means of Urgent Market Messages (UMM) in the NordPool platform. If the low-inertia event is forecast after the day-ahead market has been cleared, the downregulation is performed in the regulating power market. In these three occasions, O3 was downregulated for 166 hours and the procurement of the substitute power was performed in the regulating power market for at least 54 hours. Considering the average regulating price in summer 2018, equal to 54.06 €/MWh [19], the total cost borne by Nordic TSOs in 2018 was approx. 380 thousand Euros.

Since many more low-inertia events are to be expected in the coming years, Nordic TSOs have decided to opt for a more cost-efficient solution. Since May 2020, a new market product for frequency stability has been introduced in the Nordic System: the Fast Frequency Reserves (FFR) [20]. This product is complementary to FCR-D and does not reduce the need of these reserves. Instead, the procurement of FFR is based on inertia forecasts, and only during low-inertia periods FFR capacity is reserved. The procurement of these new reserves differs from country to country, at least until a common platform is launched. Statnett and Svenska kraftnät have implemented a seasonal procurement based on long-term forecast. In East Denmark, Energinet procures FFR on a monthly basis, with the possibility of releasing some of the reserved capacity two days before operation (D-2). Finally, Fingrid has implemented a daily procurement based on short-term forecast. The need for FFR is calculated based on an equation considering the inertia level, the dimensioning incident, the amount of FCR-D and the frequency nadir. Each TSO is responsible for procuring a fixed share of the total FFR need, respectively 14% for Energinet, 42% for Statnett, 24% for Svenska kraftnät and 20% for Fingrid [20]. Regarding summer 2020, FFR have been procured for more than 900 hours, with a total cost for Nordic TSOs equal to 3.33 million Euros.

As an alternative to these measures, HVDC lines could provide frequency support in case of disturbance. This remedial action relies on the fact that HVDC converters, equipped with fast frequency controllers, can adjust the power flow in response to frequency deviations. This control mode is referred to as Emergency Power Control (EPC). Different control strategies can be used to define the response of HVDC converters. The currently implemented strategy is based on stepwise triggers: depending on the size of the power deviation and the corresponding frequency variation, a constant amount of power is injected to improve the frequency response of the system.

The aim of the cost analysis presented in this paper is to measure the costs of these remedial actions and quantify the potential cost savings in the Nordic countries. Two future scenarios for the year 2025 are considered: "full nuclear" (2025 FN), based on the current situation with nuclear power plants fully dispatched, and "half nuclear" (2025 HN), where half of the nuclear production is replaced by wind, solar and HVDC imports. The inertia estimates come from a large number of market simulations performed by Svenska kraftnät. In order to capture a greater variety of possible weather conditions, the simulations have been performed using the meteorological conditions of the years between 1980 and 2012. The obtained results have been condensed into three cases (low, medium and high inertia) corresponding to the "worst", "most probable" and "best" cases. More information can be found in [21].

The required DI reduction, FFR and injected power by HVDC EPC, respectively, to maintain the N-1 security criterion in the Nordic system with low inertia are shown in Figure 2. These values are obtained through dynamic simulations using a single machine equivalent model of the Nordic power system. The parameters of the model are tuned based on the frequency response of the system to real disturbances that occurred in the Nordic system. The model has been validated and is currently used by Nordic TSOs in their control rooms. According to TSOs practice, the dimensioning incident is reduced by blocks of 50 MW and prolonged for few hours before and after the low-inertia event. Finally, in order to avoid the activation of FFR or HVDC EPC for small frequency deviations, a certain frequency dead band is introduced (this varies between 49.7 and 49.5 Hz). Because of this delay in the response of FFR and HVDC converters, these two measures result in being less effective than the DI reduction.

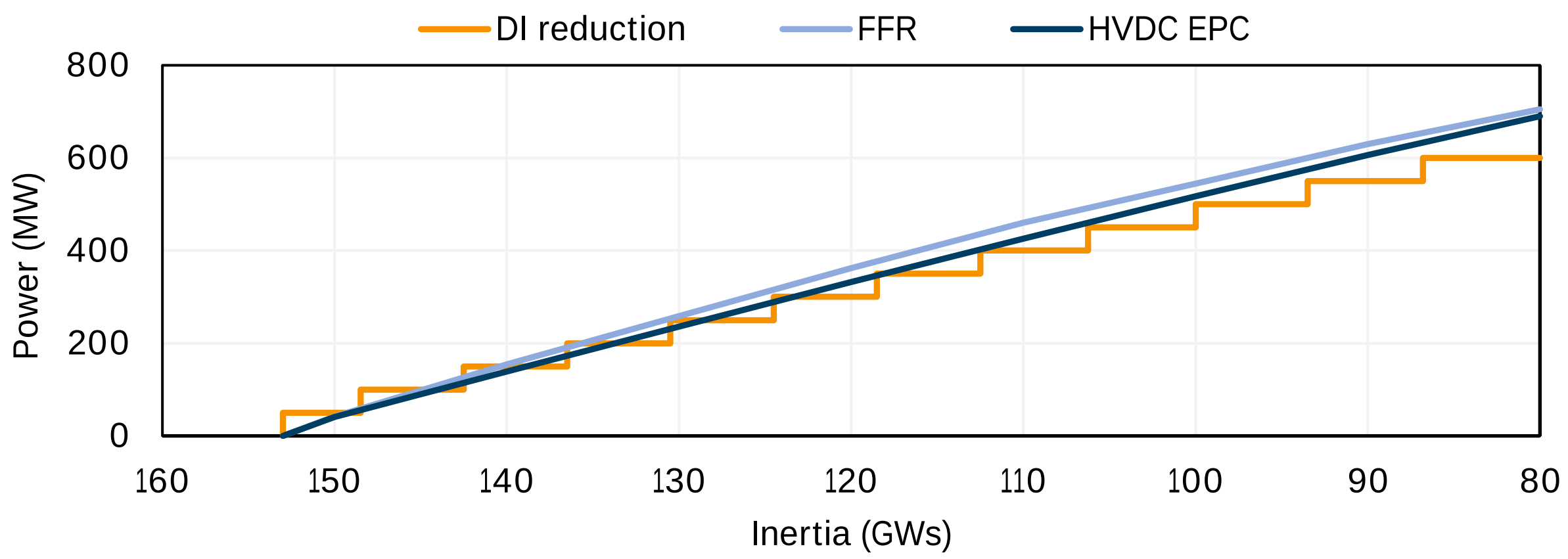

Figure 2: Required DI reduction, FFR and injected power by HVDC EPC, respectively, to maintain the N-1 security criterion in the Nordic system with low inertia.

The utilization of HVDC lines for frequency support relies on the fact that there is enough transmission capacity available on the HVDC interconnectors, and that there is a certain availability of frequency reserves on the neighboring systems. As for now, there seems to be no regulation about how these can be procured on a market basis among asynchronous areas. In order to assess what could be the cost of this remedial action, we envision a possible future situation where there is a European market for reserves, and Nordic TSOs are requested to procure the necessary primary reserves through this platform. In fact, this seems to be the direction that European countries are taking, as described in [22] for automatic activated FRR. The costs of primary frequency reserves are calculated based on historical data from the past 5 years, considering the distribution of the average price obtained via data re-sampling. Moreover, the reservation of HVDC capacity is assumed to come with a cost, as considered also in [9]. For the reservation costs, we consider that the price for reserving HVDC capacity in a specific hour is equal to the average congestion rent (the distribution is obtained via data re-sampling).

For the calculation of the costs associated with DI reduction, the downregulation of O3 is assumed to happen after the day-ahead market has been cleared and the substitute power is procured in the regulating power market for the first 24 hours of each event; OKG receives a fixed compensation per event (for the reduced efficiency) and a compensation proportional to the power limitation (for the opportunity costs). For the procurement costs of FFR, the average price of summer 2020 is used, equal to 48.95 €/MW/h. This price is calculated as the average price in the four countries (East Denmark, Norway, Sweden and Finland) weighted by the share of FFR of each country.

Table 1 presents the hours when the kinetic energy is below the required level and the corresponding remedial actions. The first thing to be noticed is that the occurrence of low inertia periods is highly impacted by the weather. In terms of hours when the kinetic energy is expected to be below the requirements, this number is almost zero in case of wet hydrological years (high inertia scenarios). On the contrary, these numbers are tripled in case of dry years (compared to an average hydrological year). This happens because during dry periods the production of hydro power plants (that are synchronous machines) is replaced by imports via HVDC (which do not contribute to the kinetic energy of the system). The same trend is observed with the size of remedial actions ("Energy" column, expressed in GWh), as they are proportional to the number of hours. Conversely, the number of events does not follow the same trend: technical limitations play an important role on the length of DI reduction and, with high frequency of occurrence of low inertia periods, the limitation on O3 is often prolonged resulting in less but longer events. Lastly, the length of each event differs depending on the remedial action. For instance, HVDC EPC and FFR are needed only for those specific hours when the kinetic energy is below the requirements, while this is not the case for DI reduction.

| | | **2025 Full Nuclear** | | | **2025 Half Nuclear** | | |
|---|---|---|---|---|---|---|---|
| | Action | Occasions | Hours | Energy | Occasions | Hours | Energy |
| **High inertia** | DI reduction | - | - | - | 1 | 60 | 5 |
| | FFR | - | - | - | - | 12 | 1 |
| | HVDC EPC | - | - | - | - | 12 | 1 |
| **Medium inertia** | DI reduction | 15 | 345 | 48 | 18 | 649 | 132 |
| | FFR | - | 165 | 11 | - | 165 | 31 |
| | HVDC EPC | - | 165 | 11 | - | 166 | 32 |
| **Low inertia** | DI reduction | 16 | 1020 | 182 | 14 | 1638 | 437 |
| | FFR | - | 411 | 28 | - | 411 | 72 |
| | HVDC EPC | - | 411 | 29 | - | 411 | 73 |

Table 1: Low inertia events and remedial actions (energy in GWh).

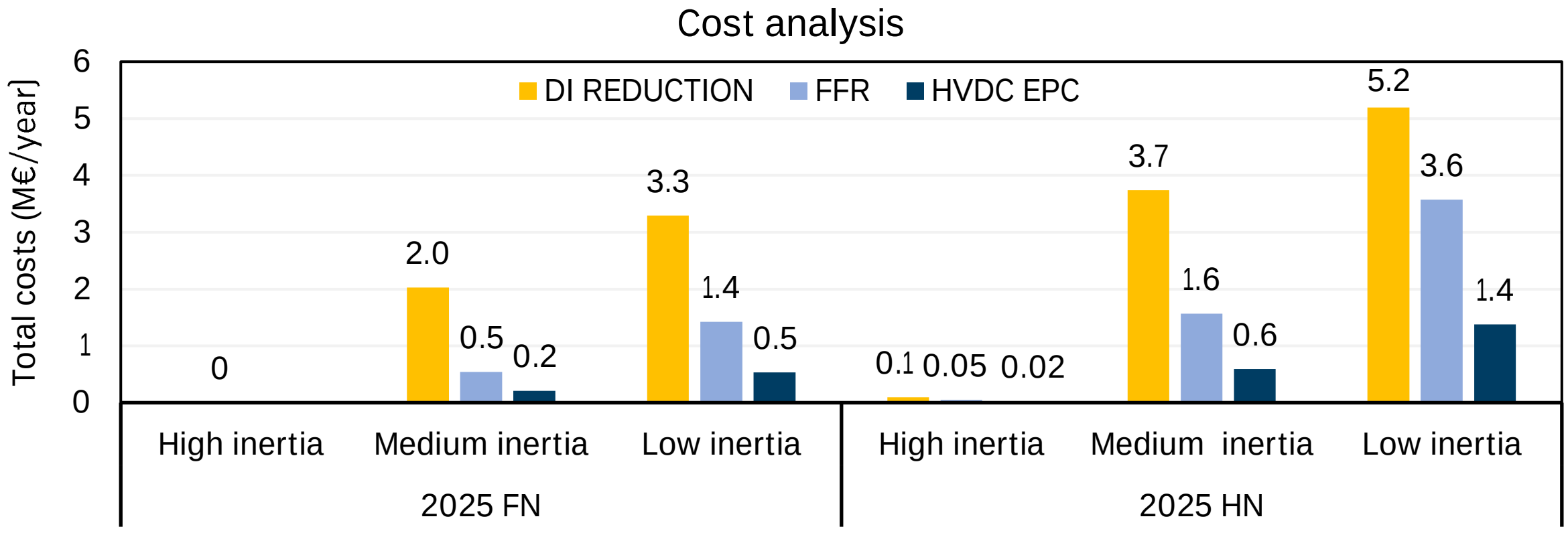

Figure 3: Total costs of the different remedial actions.

The total costs of the different remedial actions are provided in Figure 3. As expected, the DI reduction is the most expensive among the considered actions. This clearly motivates the change of paradigm from the first events in 2018. Still, however, the utilization of HVDC lines for frequency containment is the most cost-efficient solution. If HVDC EPC is implemented within the next 5 years, significant cost savings can be achieved (70-90% compared to the DI reduction case, or 60-65% compared to the FFR case).

## IV. HVDC LOSS FACTORS IN MARKET CLEARING

Nordic TSOs have proposed to introduce loss factors for HVDC lines to avoid HVDC flows between zones with zero price difference. The proposal has already gone through the first stages of the process and it is currently under investigation for real implementation in the market clearing algorithm. In [23], we developed a rigorous framework to assess this proposal; the results showed that the benefits of such a measure depend on the topology of the investigated system. In this paper, we present the results of our analyses on a detailed market model of the Nordic countries.

The focus of the analysis is on the differences between linear and piecewise-linear loss factors and between HVDC and AC+HVDC loss factors. Implementing such measures in real systems is possible: for instance, piecewise-linear loss factors are already used in real power exchanges, e.g. New Zealand Exchange (NXZ) [24], and several power markets in the US already use sensitivity factors to determine AC losses [25], [26]. Four simulations are run considering different loss factors at a time:

1. No loss factors (reference case);
2. Linear HVDC loss factors;
3. Piecewise-linear HVDC loss factors;
4. Piecewise-linear AC and HVDC loss factors.

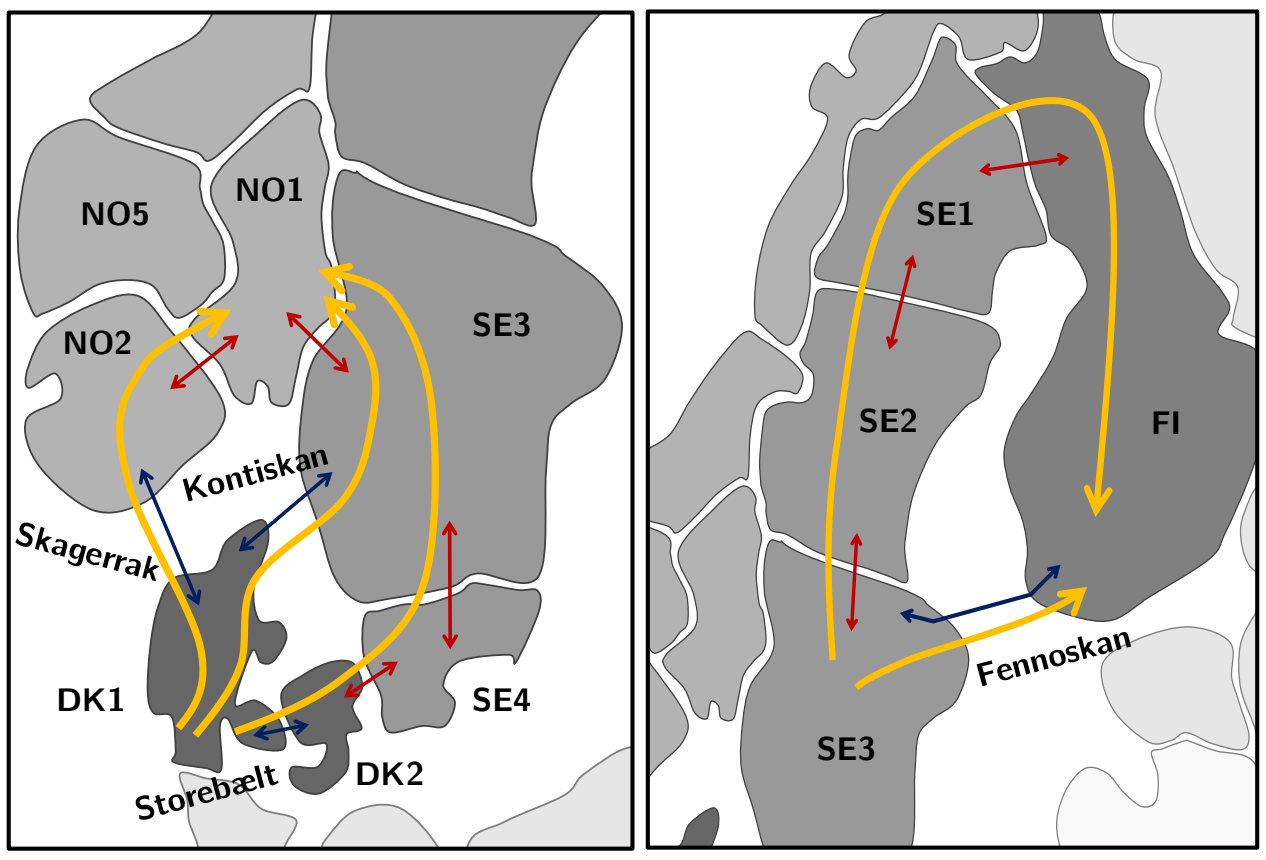

Figure 4: Example of flows on parallel HVDC paths (left) and on parallel AC and HVDC paths (right).

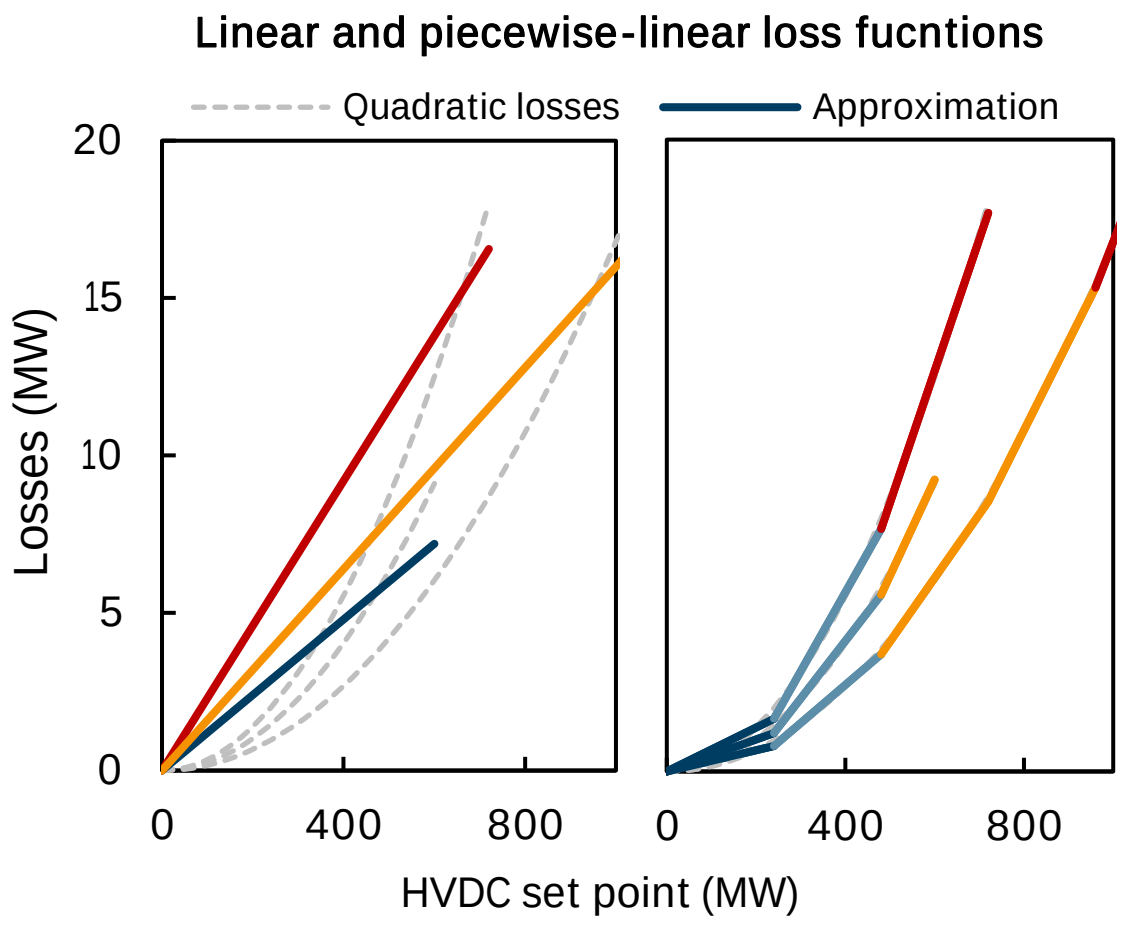

Figure 5: Linear (left) and piecewise-linear (right) loss functions for Skagerrak, Storebælt and Kontiskan. The different colors refer to the slope of the lines.

In each simulation, the market is cleared for each hour of the year (8760 instances) using data from 2017. It is important to mention that all the cost-benefit analyses are limited to the introduction of loss factors in the intra-Nordic interconnectors, that means Fennoskan, Skagerrak, Storebælt, Kontiskan and only the AC interconnectors of RG Nordic. Indeed, the power exchanges with neighboring countries are fixed to the real exchanges, and so are the flows on the interconnectors (becoming unresponsive to any change introduced by loss factors). Moreover, because of the zonal-pricing scheme, intra-zonal losses are not considered in the analysis: all the presented results are limited to losses on the interconnectors.

With the inclusion of HVDC loss factors in the market, HVDC losses are implicitly considered when the market is cleared. Since losses appear in the power balance equation, they represent an extra cost and the solver will try to minimize them. Given that only HVDC losses are considered, the solver will use HVDC interconnectors only if necessary, i.e. in case of congestions in the AC system or for exchanges between asynchronous regions. For the same reason, when forced to use HVDC interconnectors, the solver will look at which path produces the least amount of losses. In case of linear loss factors, the slope of the linear loss functions is the discriminating factor. This might become a problem in a situation with different parallel HVDC paths, as it is the case, for example, of Skagerrak, Kontiskan and Storebælt in Western Denmark (Figure 4 - left). In such a situation, the solver will direct the flow over the line with the smallest slope (in the left chart of Figure 5, the blue one) and only when its capacity is fully utilized it will start directing the flow towards the line with the second smallest slope (the orange one), and finally towards the remaining line (the red one).

With piecewise-linear loss functions, the solver finds the path that produces the least amount of losses by moving back and forth from one loss function to the other. As with linear loss factors, it will start with the HVDC line with the smallest slope. However, since the slope changes in the next segment, the solver will start directing the power flow towards other lines if the slope of those segments is smaller (in the right chart of Figure 5, all the blue segments). It will move back to the first line only when there are no other segments with smaller slopes, i.e. it will move to orange segments when there are no more blue segments, and so on. In this way, the quadratic nature of losses is better represented, allowing the solver to identify the best path and better distribute the power flows among the HVDC lines.

Similarly, with the inclusion of only HVDC loss factors, the solver will see HVDC lines as expensive alternatives to AC lines, whose losses are not considered when the market is cleared: if there exist parallel AC and HVDC paths, the solver will always prefer the AC option. This is the case, for example, of Fennoskan, the HVDC link connecting Sweden and Finland (Figure 4 - right). In this case, if implicit grid loss is implemented on Fennoskan and not on the AC interconnectors SE3-SE2, SE2-SE1 and SE1-FI, the solver will always try to reroute the power across the AC path. However, losses are produced in the AC system as well and, by reducing the flow on some HVDC interconnectors, we might disproportionally increase losses in the AC system. The only way to minimize losses and maximize social benefits is to include loss factors for AC interconnectors as well. By doing so, the solver will be able to identify the path producing the least amount of losses.

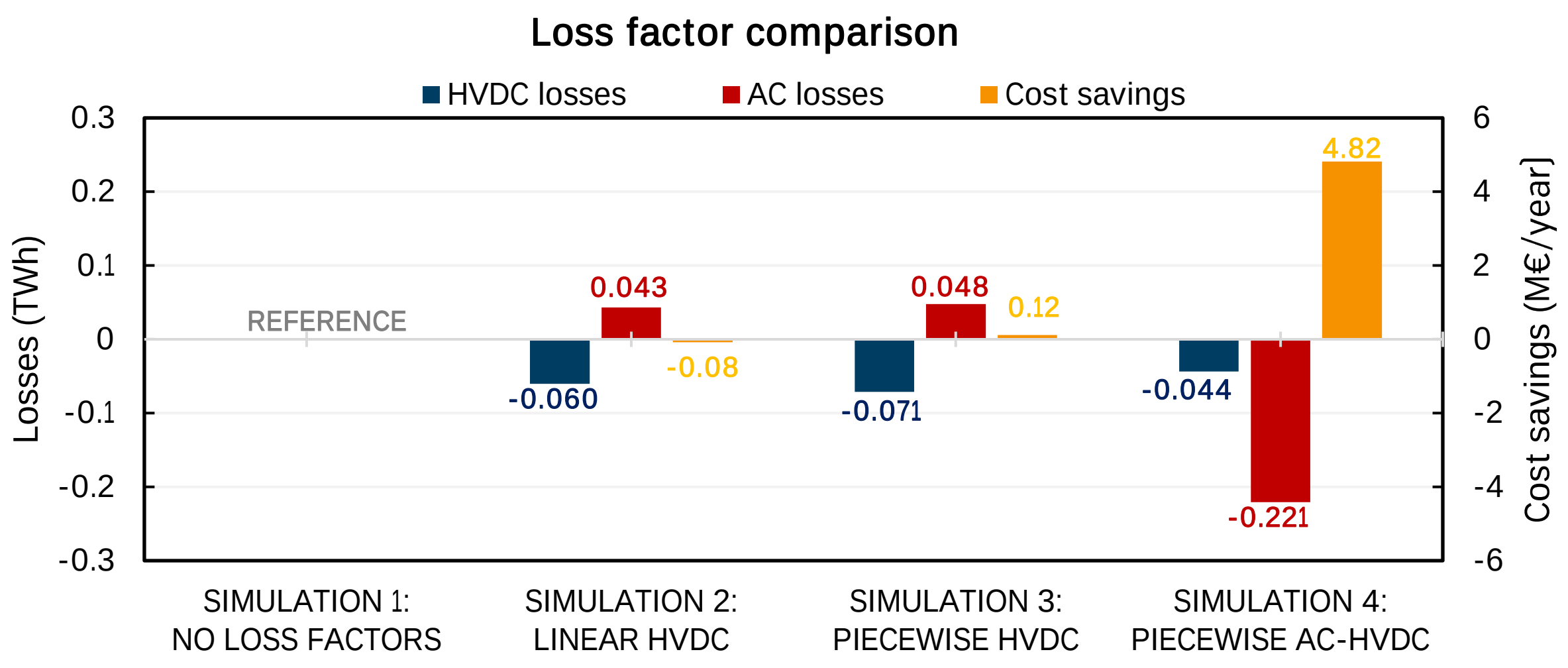

Figure 6: Comparison of the simulation with focus on AC and HVDC losses. The reference case is the simulation without loss factors.

The comparison of the four simulations is shown in Figure 6, where blue bars represent HVDC losses, red bars AC losses and yellow bars cost savings. As expected, in simulation 2 and 3, the reduction of HVDC losses comes together with an increase of AC losses. The net reduction of losses is positive, meaning that the introduction of only HVDC loss factors can be beneficial; however, the resulting cost savings in simulation 2 are negative. This happens because linear loss factors result in a bad approximation of losses which are often overestimated, meaning that unnecessary power is provided by generators (at a higher cost for society). This does not happen with piecewise-linear loss factors because they better represent HVDC loss functions.

The results of simulation 4 show that it is possible to decrease the sum of AC and HVDC losses by 12% (compared to simulation 1, where losses on the interconnectors amount to 2.42 TWh) by introducing piecewise-linear loss factors for AC and HVDC interconnectors, while this is limited to 0.7% with only linear HVDC loss factors and to 0.9% with only piecewise-linear HVDC loss factors. Concerning the cost savings, they increase moving from left to right in Figure 6, showing the progressive benefit of having piecewise-linear loss factors and AC loss factors. In particular, simulation 4 with piecewise-linear loss factors for both AC and HVDC interconnectors results in cost savings of 4.82 million euros per year.

## V. CONCLUSION

In the Nordic countries, more than 10 interconnectors are formed by HVDC links, and many new projects are under construction or under investigation. Based on this consideration, this paper explores the potential benefit of using HVDC links for the exchange of ancillary services between countries and investigates different solutions for the inclusion of HVDC losses in the market.

The first analysis was motivated by the low-inertia events occurred in 2018, during which Svenska kraftnät had to reduce the output of Oskarshamn 3 to guarantee N-1 security.

Because of the increasing frequency of occurrence of low inertia events, Nordic TSOs have introduced Fast Frequency Reserves in summer 2020. The performance of this new product has been verified in this paper for the studied incident. However, our results show that, if HVDC is used in the form of Emergency Power Control, these costs could be more than halved.

The second analysis comes from the proposition of Nordic TSOs of including linear loss factors for HVDC lines to avoid flows between zones with zero price difference. Our results show that there is room for improvement in two directions. First, by using piece-wise linear loss factors. This would lead to a better representation of loss functions and to an optimal distribution of power flows, resulting in a further decrease of losses and higher cost savings. Second, by also introducing AC loss factors. This would allow for the identification of the optimal paths that leads to the least amount of losses, maximizing cost savings.

## VI. ACKNOWLEDGMENT

This work was supported by the multiDC project, funded by Innovation Fund Denmark under Grant 6154-00020B.